\documentclass[11pt,twoside]{article}
\usepackage{latexsym}
\usepackage{amssymb,amsbsy,amsmath,amsfonts,amssymb,amscd,color}
\usepackage{epsfig, graphicx}
\usepackage{hyperref}
\usepackage{dsfont}
\newcommand{\commentout}[1]{}
\newcommand{\R}{\mathbb{R}}

\newcommand {\eps}  {\varepsilon}

\newcommand {\Chi} {{\bf \raise 2pt \hbox{$\chi$}} }

\newcommand {\dv}  { {\rm div} }
\newcommand {\f}   {\frac}
\newcommand {\p}   {\partial}
\newcommand{\dis}{\displaystyle}
\newcommand{\beq}{\begin{equation}}
\newcommand{\eeq}{\end{equation}}
\newcommand{\bea} {\begin{array}{rl}}
\newcommand{\eea} {\end{array}}
\newcommand{\bepa}{\left\{ \begin{array}{l}}
\newcommand{\eepa} {\end{array}\right.}
\newtheorem{theorem}{Theorem}[section]
\newtheorem{lemma}[theorem]{Lemma}

\newcommand{\qed}{\hfill$\square$\vspace{0.3cm}}
\title{\Large \bf Derivation of a Hele-Shaw type system from a cell model with active motion}

\author{
Beno\^ \i t Perthame\thanks{UPMC Univ Paris 06 and CNRS UMR 7598, Laboratoire Jacques-Louis Lions, F-75005, Paris, France. Email B. P.: benoit.perthame@ljll.math.upmc.fr, Email N. V.: vauchelet@ann.jussieu.fr}  \thanks{INRIA Paris Rocquencourt EPC BANG. France.}
\and   Fernando Quir\'os\thanks{Departamento de Matem\'aticas, Universidad
Aut\'onoma de Madrid,  28049-Madrid, Spain. Supported by  Spanish project
MTM2011-24696. Email:
fernando.quiros@uam.es}
\and  Min Tang\thanks{Department of mathematics, Institute of Natural Sciences and MOE-LSC. Shanghai Jiao Tong University, China. Email: tangmin1002@gmail.com}
\and Nicolas Vauchelet\footnotemark[1] \footnotemark[2]
}

\date{\today}

\begin{document}
\maketitle
\pagestyle{plain}
\begin{abstract}
We formulate a Hele-Shaw type  free boundary problem for  a tumor growing under the combined effects of
pressure forces, cell multiplication and active motion, the latter being the novelty of the present paper.  This new ingredient is considered here as  a standard diffusion process.
The free boundary model is derived from a description at the cell level using the asymptotic of a stiff pressure limit.

Compared to the case when active motion is neglected, the pressure satisfies the same complementarity Hele-Shaw type formula. However, the cell density is smoother (Lipschitz continuous),  while   there is a deep  change in the  free boundary velocity, which is no longer given by the gradient of the pressure, because some kind of \lq mushy region' prepares  the tumor invasion.

\end{abstract}


\pagenumbering{arabic}

\section{Introduction}
\setcounter{equation}{0}

Among the several models now available to deal with cancer development, there is a class, initiated in the 70's by Greenspan \cite{greenspan}, that considers that cancerous cells multiplication is limited by nutrients
(glucosis, oxygen) brought by blood vessels. Models of this class rely on two kinds of descriptions; either they describe the dynamics of cell
population density \cite{byrne-chaplain} or they consider the \lq
geometric' motion of the tumor through a free boundary problem; see
\cite{cui_escher, cui, friedman_hu, Lowengrub_survey} and the references
therein. In the latter kind of models the stability or  instability of the free boundary is an important issue that  has attracted attention, \cite{ciarletta,friedman_hu}.

The first stage,  where growth is limited by nutrients,  lasts until the tumor reaches the size of
$\approx 1$mm; then, lack of food leads to cell necrosis which
triggers neovasculatures development \cite{chaplain} that supply the
tumor with enough nourishment. This has motivated a new generation of
models where growth is limited by the competition for space
\cite{Bru}, turning the modeling effort towards
mechanical concepts, considering tissues as multiphasic fluids (the
phases could be intersticial water, healthy and tumor cells,
extra-cellular matrix \dots)  \cite{byrne-drasdo, byrne-preziosi,
preziosi_tosin,Bellomo1, RCM}. This point of view is now sustained by experimental
evidence \cite{JJP}. The term \lq homeostatic pressure', coined
recently, denotes the lower pressure that prevents cell
multiplication by contact inhibition.

In a recent paper \cite{PQV} the authors explain how asymptotic analysis can link
the two main  approaches, cell density models and free boundary
models, in the context of fluid mechanics for the
simplest cell population density model, proposed
in~\cite{byrne-drasdo}, in which the cell population density evolves under pressure forces and cell
multiplication. The principle of the derivation is to use the stiff limit in the pressure law of state, as treated in several papers; see for instance  \cite{BI2, GQ2} and the references therein.

Besides mechanical motion induced by pressure, for some types of cancer cells it is  important to take into account active motion; see \cite{BOBAM, Drasdo_H, SLCF}.  In the present paper we extend the asymptotic analysis of \cite{PQV} to a model that includes such an ingredient. We examine the specific form of the Hele-Shaw limit and draw qualitative conclusions on the  behaviour of the solutions in terms of regularity and free boundary velocity.

\section{Notations and main result}
\setcounter{equation}{0}

Our model of tumor growth incorporates active motion of cells thanks to a diffusion term,
\begin{equation}
\dis \partial_t n_k - \dv \big( n_k \nabla p_k \big) - \nu \Delta n_k  = n_k G\big(p_k \big),  \qquad (x,t)\in Q:= \R^d\times(0,\infty).
\label{eq:n}
\end{equation}
The variable $n_k$ represents the density of tumor cells, and the variable $p_k$ the pressure, which is considered to be given by a homogeneous law (written with a specific coefficient so as to simplify notations later on)
\begin{equation}
\dis p_k(n) = \frac{k}{k-1} n^{k-1}.
\label{eq:p}
\end{equation}
Hence, we are dealing with a porous medium type equation; see~\cite{JLV} for a general reference on such problems.
We complement this system with an initial  condition that is supposed to satisfy
\begin{equation} \label{hypInit} \left\{\begin{array}{l}
n_k(x,0)=n^{ini}(x) > 0, \qquad  n^{ini}\in L^1(\mathbb{R}^d)\cap L^\infty(\mathbb{R}^d),
\\[6pt]
p_k^{ini}:=\f{k}{k-1}(n^{ini})^{k-1} \leq P_M .
\end{array} \right.\end{equation}

In a purely mechanical view, the pressure-limited growth is
described by the function $G$, which satisfies
\beq\label{hypG}
G'(\cdot ) < 0 \qquad \mbox{and }\quad G(P_M)=0,
\eeq
for some $P_M>0$, usually called the homeostatic pressure; see \cite{byrne-drasdo,JJP}.

Many authors use another type of models, namely  free boundary problems on the tumor region $\Omega(t)$. Our purpose is to make a rigorous derivation of one of such  models from \eqref{eq:n}, \eqref{eq:p}. As it is wellkown, for $\nu =0$ this is possible in  the asymptotics $k$ large. This is connected, in fluid mechanics, to the Hele-Shaw equations; a complete proof of the derivation is provided in \cite{PQV}. Typically the limit of the cell density is an indicator function for each time $t>0$,  $n_\infty = \mathds{1}_{\Omega(t)}$, if this is initially true, and the problem is reduced to describing the velocity of the boundary $\p \Omega(t)$.

Our aim is thus to understand what is the effect of including active motion, that is,  $\nu > 0$. We will show that both the density and the pressures have limits, $n_\infty$ and $p_\infty$, as $k\to\infty$ that satisfy
\beq
\partial_t n_\infty - \dv \big( n_\infty \nabla p_\infty \big) - \nu \Delta n_\infty  = n _\infty G\big(p_\infty \big).
\label{eq:ninfty2}
\eeq
Compared with the case $\nu=0$ considered in \cite{PQV}, a  first major difference is that now the cell density $n_k$ is smooth, since equation~\eqref{eq:n} is non-degenerate when $\nu>0$. Is that translated into more regularity for the limit density? We will show that this is indeed the case. Though the limit density satisfies
$$
0\le n_\infty \leq 1,
$$
it is not an indicator function any more, and its time derivate $\partial_t n_\infty$ is a function, while it is only a measure when $\nu=0$.
As for the pressure, we will establish that we still have
$$
n_\infty=1  \hbox{ in }   \Omega(t)=\{p_\infty(t) >0\},
$$
or in other words  $p_\infty\in P_\infty(n_\infty)$, with $P_\infty$ the limiting monotone
graph
\beq
P_\infty(n) = \left\{\begin{array}{ll}
0, & 0\leq n < 1, \\[2mm]
[0,\infty),\qquad & n=1.
\end{array}
\right.
\label{eq:graph}
\eeq
Furthermore,  multiplying equation \eqref{eq:n} by $p_k'(n_k)$ leads to
$$
\partial_t p_k- n_k p_k'(n_k) \Delta p_k - |\nabla p_k|^2 - \nu \Delta p_k = n_k p_k'(n_k)  G\big(p_k \big) - \nu p_k''(n_k) |\nabla n_k|^2,
$$
and for the special case $p_k= \frac{k}{k-1} n_k^{k-1}$ at hand we find
\begin{equation}    \label{eq:pk}
\partial_t p_k- (k-1) p_k \Delta p_k - |\nabla p_k|^2 - \nu \Delta p_k =
(k-1) p_k G\big(p_k \big) - \nu \frac{(k-2) \nabla p_k\cdot \nabla n_k}{n_k}.
\end{equation}
Therefore, the \lq complementary relation'
\beq
- p_\infty \Delta p_\infty = p_\infty G\big(p_\infty\big) - \nu \frac{\nabla p_\infty \cdot\nabla n_\infty}{n_\infty},
\label{eqp1}
\eeq
is expected in the limit. However,  $\nabla p_\infty$ vanishes unless $p_\infty>0$, in which case $n_\infty =1$, therefore $\nabla n_\infty =0$.
Thus, the equation on $p_\infty$ ignores the additional term coming from active motion and reduces to the same Hele-Shaw equation for the  pressure that holds when $\nu=0$, namely
\beq
 p_\infty \big[ \Delta p_\infty + G\big(p_\infty \big) \big] =0.
\label{eqp2}
\eeq
Let us remark that this does not mean that active motion has no effect in the limit. Though the pressure equation is the same one as for the case $\nu=0$, the  free boundary $\p \Omega(t)$ is not expected to move with the usual Hele-Shaw rule $V= -\nabla p_\infty$, but with a faster one; see Section~\ref{sec:fb.velocity} for a discussion on the speed of the free boundary.

The above heuristic discussion can be made rigorous.
\begin{theorem}\label{th1}
Let $T>0$ and $Q_T=\mathbb{R}^d\times (0,T)$.
Assume \eqref{hypInit}, \eqref{hypG} and that the initial data satisfies
$\p_t n^{ini}\geq 0$.
Consider  a weak solution $(n_k,p_k)$  of \eqref{eq:n}--\eqref{eq:p}. Up to extraction of a subsequence,  $(n_k,p_k)_k$ converges strongly in $L^p(Q_T)$, $1\leq p < \infty$, to limits
$$
n_\infty\in C\big([0,\infty); L^1(\mathbb{R}^d)\big)\cap
L^\infty((0,T);H^1(\mathbb{R}^d)),\qquad p_\infty\in L^\infty((0,T);H^1(\mathbb{R}^d)),
$$
such that
$0\leq n_\infty \leq  1$, $n_\infty(0)=n^{ini}$, $0\leq p_\infty \leq P_M$, $p_\infty\in P_\infty(n_\infty)$, where $P_\infty$ is the Hele-Shaw monotone graph given in \eqref{eq:graph}. Moreover, the pair $(n_\infty,p_\infty)$ satisfies
on the one hand \eqref{eq:ninfty2}, and on the other hand the Hele-Shaw type equation
\beq
\partial_t n_\infty - \Delta p_\infty -\nu \Delta n_\infty  = n _\infty G\big(p_\infty \big),
\label{eq:ninfty}
\eeq
and the complementarity relation~\eqref{eqp2} for  almost every $t>0$, all three equations
in the weak sense. The time derivatives of the limit functions satisfy
$$
\p_t n_\infty,\;\p_t p_\infty \in {\cal M}^1(Q_T),\qquad \partial_t n_\infty,\;\partial_t p_\infty\ge0.
$$
\end{theorem}

\begin{figure}
\begin{center}
\includegraphics[width=0.25\textwidth]{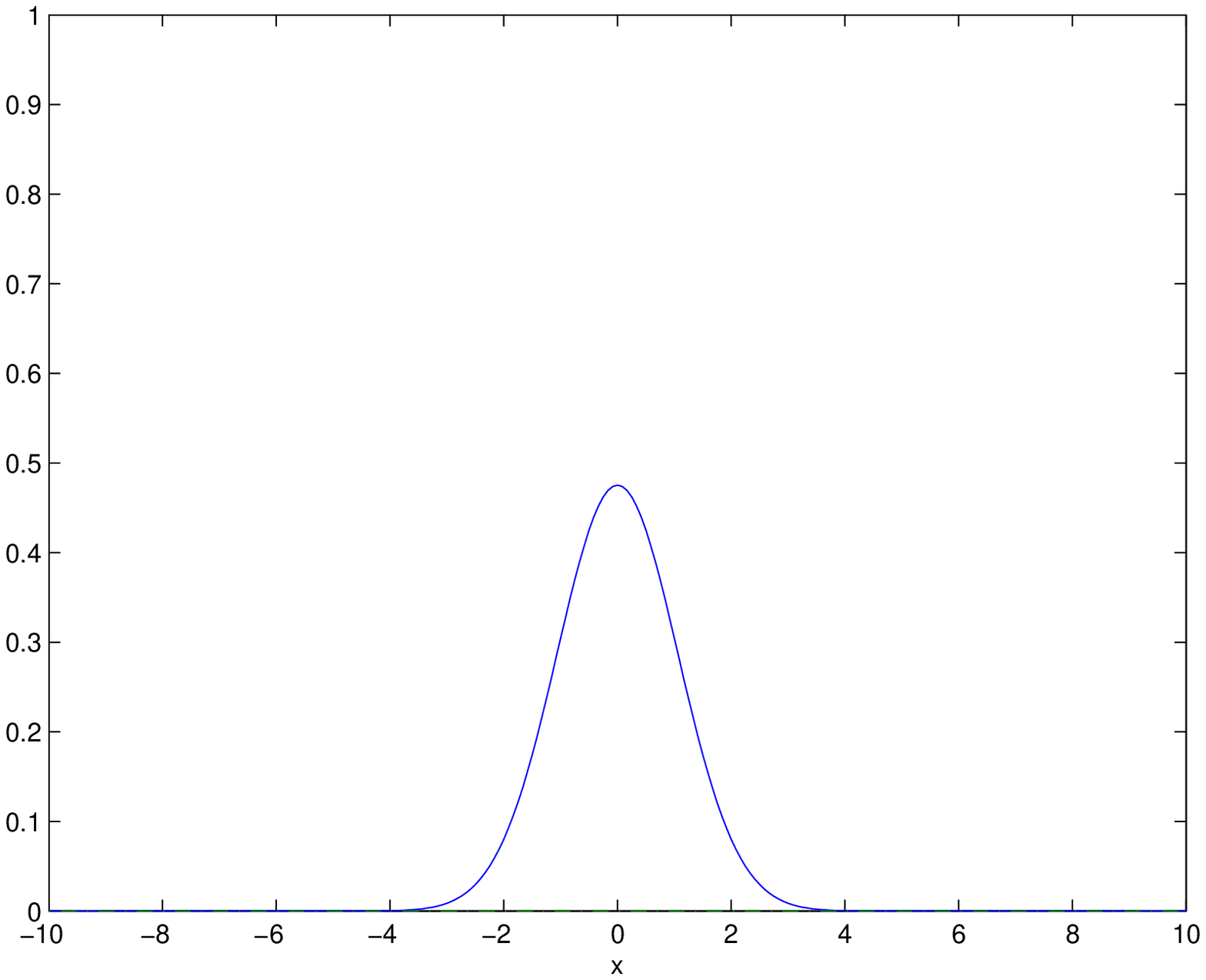} \hspace{-5mm}
\includegraphics[width=0.25\textwidth]{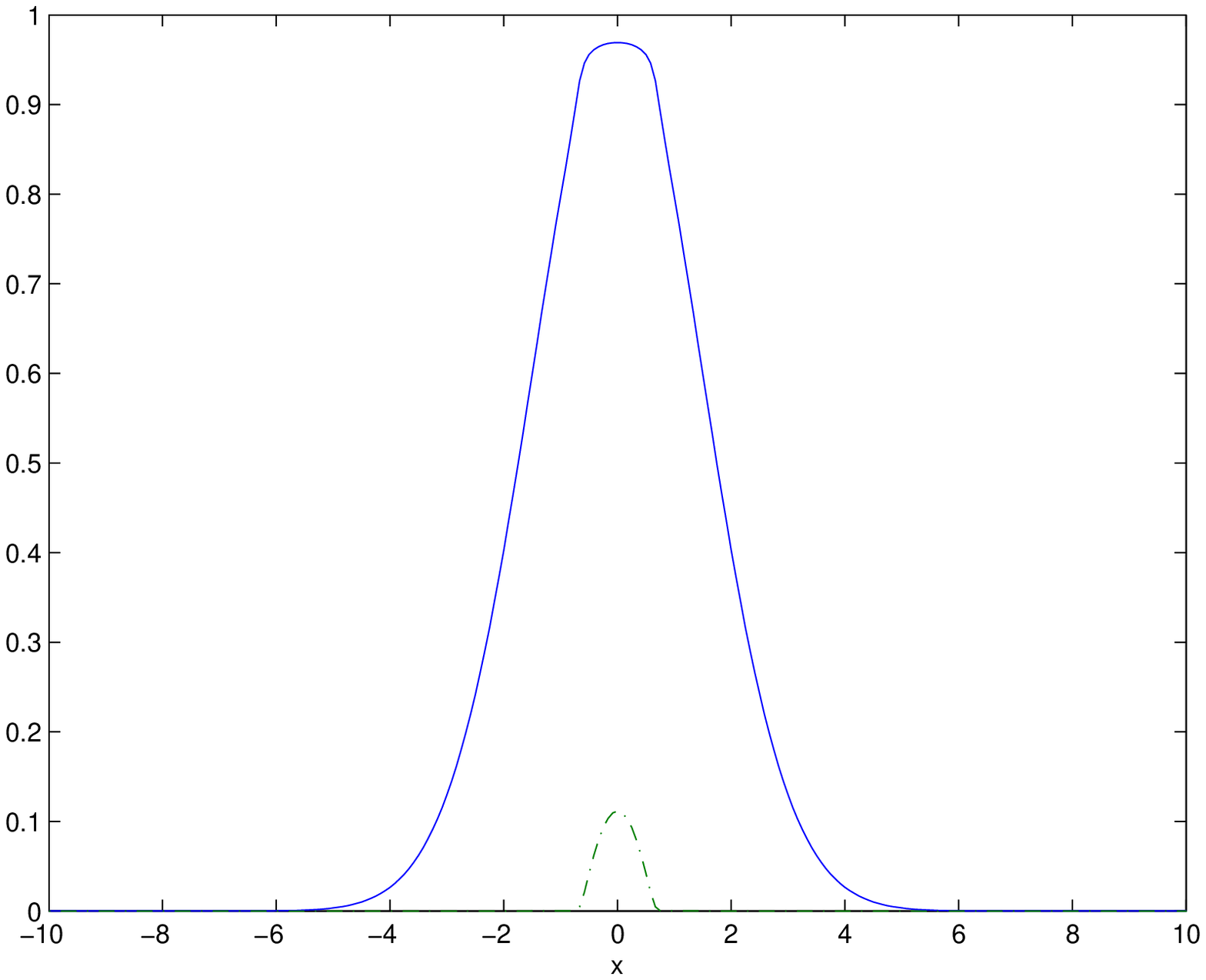}\hspace{-5mm}
\includegraphics[width=0.25\textwidth]{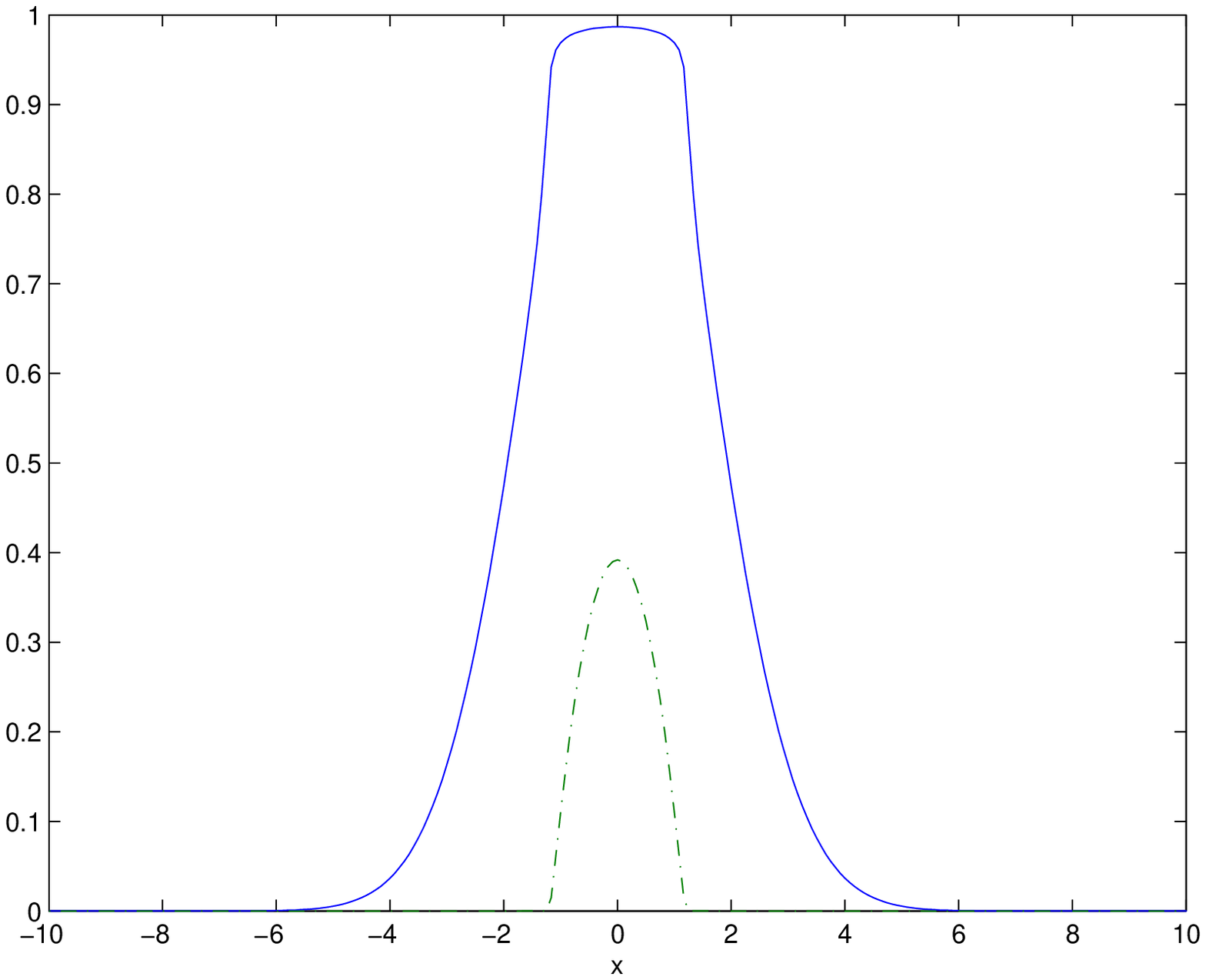}\hspace{-5mm}
\includegraphics[width=0.25\textwidth]{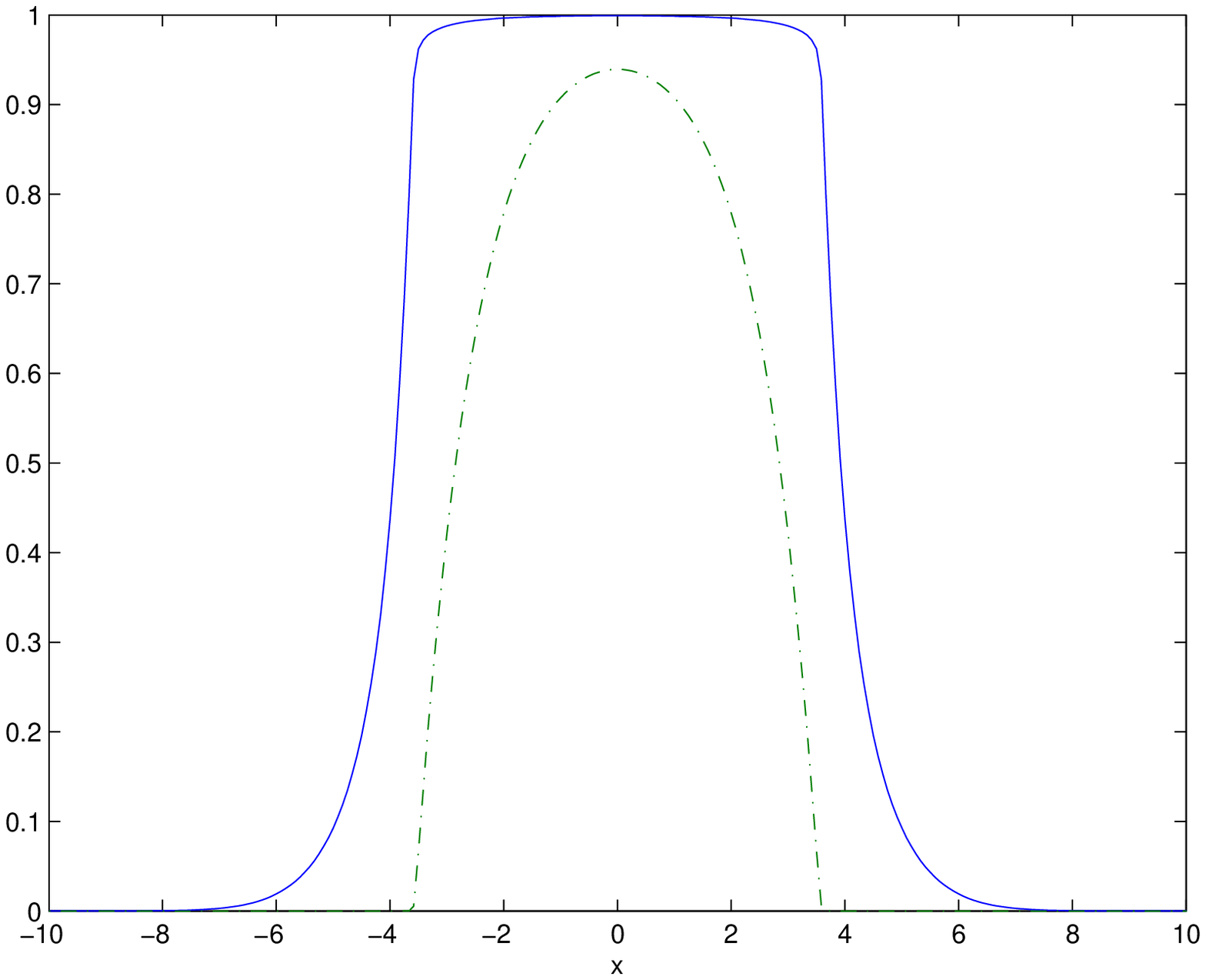}
\vspace{-5mm}
\caption{First steps of the initiation of the free boundary. Results obtained
thanks to a discretization of the system \eqref{eq:n}--\eqref{eq:p}  with $k=100$ and $\nu=0.5$.
The density $n$ is plotted in solid line whereas the pressure $p$ is represented
in dashed line. The pressure $p$ has the same shape as in the classical
Hele-Shaw system with growth. However the density $n$ is smoother.
}\label{fig:act}
\end{center}
\end{figure}

To illustrate this behaviour, we present numerical results obtained
thanks to a discretization with finite volume of system
\eqref{eq:n}--\eqref{eq:p} in the case $k=100$, $\nu=0.5$ and with $G(p)=1-p$.
We display in Figure \ref{fig:act} the first steps of the formation
of a tumor which is initially given by a small bump. As expected, we notice
that the density $n$ is smooth. The shape of the pressure $p$ at the
place where $n=1$ is similar to the one observed for the classical
Hele-Shaw system (see e.g. \cite{PQV}).

The rest of the paper is organized as follows. We begin in Section~\ref{sec:estimates} with some uniform (in $k$)  a priori estimates which are necessary for strong compactness. Then, in Section \ref{sec:proof} we prove the main statements in Theorem \ref{th1}. The most delicate part, establishing \eqref{eqp2}, is postponed to Section~\ref{sec:hs}. After proving uniqueness for the limit problem in Section~\ref{sec:uniq}, we end with a final section devoted to discuss further regularity issues and speed of the boundary of the tumor zone.

\section{Estimates}
\label{sec:estimates}
\setcounter{equation}{0}

To begin with, we gather   in the following statement all the a priori  estimates that we need later on.

\begin{lemma}\label{estim} With the assumptions and notations in Theorem \ref{th1}, the weak solution $(n_k,p_k)$ of \eqref{eq:n}--\eqref{eq:p} satisfies
$$
0\le n_k   \leq  \Big(\f{k-1}{k} P_M \Big)^{1/(k-1)} \underset{k \to \infty}{\longrightarrow} 1,
 \qquad 0\le p_k \leq P_M,
$$
$$
\int_{\mathbb{R}^d} n_k(t)\leq e^{G(0) t } \int_{\mathbb{R}^d}  n^{ini}, \qquad \int_{\mathbb{R}^d} p_k(t)\leq C e^{G(0) t } \int_{\mathbb{R}^d}  n^{ini}.
$$
with $C$ a constant independent of $k$.
Furthermore, there exists a uniform (with respect to $k$) nonnegative constant such that
\begin{equation}
  \label{eq:bornpn}
  \int_{\mathbb{R}^d}\left(\nu |\nabla n_k|^2+ k n_k^{k-1}|\nabla n_k|^2+ |\nabla p_k|^2\right)(t)\leq C\left(T, \|n^{ini}\|_{L^1(\mathbb{R}^d)\cap L^\infty(\R^d)} \right)\quad \text{for all }t\in (0,T).
\end{equation}
Finally,
\begin{equation*}
  \label{eq:sign}
\p_t n_k,\p_t p_k\geq 0, \quad \p_t n_k \;  \hbox {is bounded in}\;  L^\infty((0,T);L^1(\R^d)), \quad  \; \p_t p_k  \;  \hbox {is bounded in}\;  L^1(Q_T).
\end{equation*}
\end{lemma}

\noindent\emph{Proof. } \textit{Estimates on $n_k$ and $p_k$. } The $L^\infty(Q)$ bounds are a consequence of standard comparison arguments for \eqref{eq:n} and \eqref{eq:pk}.
The  $L^\infty((0,T);L^1(\R^d))$ bound for $n_k$ can be obtained by integrating \eqref{eq:n} over $\R^d$ and then  using \eqref{hypG}. The  $L^\infty((0,T);L^1(\R^d))$ bound for $p_k$ now follows from the relation between $p_k$ and $n_k$.

\medskip

\noindent\textit{Estimates on the time derivatives. }
We introduce the quantity
\begin{equation}
\Sigma(n_k) =  n_k^k  + \nu n_k , \qquad \Sigma'(n_k)= k n_k^{k-1} + \nu.
\label{def:sigmak}
\end{equation}
The density equation \eqref{eq:n} is rewritten in terms of this new variable as
\begin{equation}
\partial_t n_k - \Delta \Sigma(n_k) = n_k G(p_k).
\label{eq:sigmak}
\end{equation}
Using the notation $\Sigma_k=\Sigma(n_k)$ and
multiplying the above equation by $\Sigma'(n_k)$, we get
\begin{equation}
\label{eq:sigma}
\partial_t \Sigma_k- \Sigma_k' \Delta \Sigma_k = n_k \Sigma_k'  G\big(p_k \big).
\end{equation}
Let $w_k= \partial_t \Sigma(n_k)$. Notice that $\mbox{sign}\,(\p_tn_k)=\mbox{sign}\,(w_k)$.
A straightforward computation yields
$$
\partial_t w_k -\Sigma'_k \Delta w_k =
\p_t n_k \Sigma_k'' \big( \Delta \Sigma_k + n_k G(p_k)\big) +
\p_t n_k \Sigma_k'G(p_k) +\p_tn_k \Sigma'_k k n_k^{k-1} G'(p_k).
$$
By using that $w_k=\Sigma'_k\p_tn_k$ and
$\Sigma'(n_k)=k n_k^{k-1}+\nu \geq \nu > 0$,
the right hand side of the above equation can be written in a more handful way as
$$
\partial_t w_k - \Sigma_k' \Delta w_k =
w_k \Big( \f{\Sigma''_k}{\Sigma'_k}
\big(\Delta \Sigma_k+n_k G(p_k) \big) +G(p_k) + k n_k^{k-1} G'(p_k) \Big).
$$
Since this equation preserves positivity and $\mbox{sign}\, (w_k(0))=\mbox{sign}\, (\p_t n_k^{ini})   \geq 0$, we conclude that $ w_k\geq 0$, that is, $ \p_t n_k  \geq 0$. The relation between $p_k$ and $n_k$ then immediately yields $\partial_t p_k\ge 0$.

Now that we know that the time derivatives have a sign, bounds for them follow easily. Indeed, using~\eqref{eq:n}, we get
$$
\|\p_t n_k(t)\|_{L^1(\R^d)} = \f{d}{dt} \int_{\R^d} n_k(t) \leq G(0) \|n_k(t)\|_{L^1(\R^d)}.
$$
This gives the bound on $\p_t n_k$ in $L^\infty([0,T];L^1(\R^d))$.
For $\p_t p_k$ we write
$$
\| \p_t p_k\|_{L^1(Q_T)} =
\int_0^T \f{d}{dt}\left( \int_{\R^d} p_k(t)\right)\, dt  \leq  \int_{\R^d} p_k(T).
$$
This last expression is uniformly bounded in $k$.

\medskip

\noindent\textit{Estimates on the gradients. }
We multiply equation \eqref{eq:n} by $n_k$, integrate over $\R^d$
and use integration by parts for the diffusion terms,
\begin{equation*}
\label{eq:estimate.gradient.n}
\int_{\R^d}(n_k\partial_t n_k)(t) + \int_{\R^d} \big(k n_k^{k-1}|\nabla n_k|^2
+ \nu|\nabla n_k|^2\big)(t) = \int_{\R^d} (n_k^2 G(p_k))(t)
\leq G(0)  \int_{\R^d}n_k^2(t).
\end{equation*}
Since both $n_k$ and $\partial_t n_k$ are nonnegative, we immediately obtain the estimate on the first two terms in \eqref{eq:bornpn}. On the other hand, integrating equation \eqref{eq:pk}, we deduce
\begin{equation*}
\label{eq:estimate.gradient.pressure}
\begin{array}{ll}
\displaystyle\int_{\R^d} \partial_t p_k(t) + (k-2)\int_{\R^d} \big(|\nabla p_k|^2
+\nu k n_k^{k-3} |\nabla n_k|^2\big)(t) &
\dis = (k-1)\int_{\R^d}( p_k G(p_k))(t) \\[3mm]
&\dis \leq (k-1) G(0) \int_{\R^d} p_k(t).
\end{array}
\end{equation*}
Since $\partial_t p_k\ge 0$, we easily obtain
the $L^2$ bound on $\nabla p_k$ in \eqref{eq:bornpn}.
\qed

\section{Proof of Theorem \ref{th1}}
\label{sec:proof}
\setcounter{equation}{0}

In this section we prove all the statements in Theorem \ref{th1} except the one concerning the complementarity relation for the pressure, equation~\eqref{eqp2}, whose proof is postponed to the next section.

\medskip

\noindent \emph{Strong convergence and bounds. }
Since the families $n_k$ and $p_k$ are bounded in $W^{1,1}_{\rm loc}(Q)$, we have strong convergence in $L^1_{\rm loc}$ both for $n_k$ and $p_k$. To pass from local convergence to convergence in $L^1(Q_T)$,
we need  to prove that the mass in an initial strip $t\in[0,1/R]$
and in the tails $|x|>R$ are uniformly (in $k$) small if $R$ is
large enough. The control on the initial strip is immediate using
our uniform, in $k$ and $t$, bounds for
$\|n_k(t)\|_{L^1(\mathbb{R}^d)}$ and
$\|p_k(t)\|_{L^1(\mathbb{R}^d)}$. The tails for the densities $n_k$ are controlled using the equation, pretty in the same way as it was done for the case $\nu=0$; see~\cite{PQV} for the details. The control on the tails of the pressures $p_k$ then follows from the relation between $p_k$ and $n_k$.
Strong convergence in $L^p(Q_T)$ for $1<p<\infty$ is now a consequence of the uniform bounds for $n_k$ and $p_k$.

Thanks to the a priori estimates proved above,
we also have that
$(\nabla n_k)_k$ and $(\nabla p_k)_k$ converge weakly in $L^2(Q_T)$, and
$$
0\leq n_\infty \leq 1,\quad  n_\infty, \;   p_\infty \in L^\infty((0,T);H^1(\R^d)),\quad
\partial_t n_\infty,\; \partial_t p_\infty\in\mathcal{M}^1(Q_T), \quad   
 \p_t n_\infty,\; \p_t p_\infty  \geq 0.
$$

\medskip

\noindent\emph{Identification of the limit. }
To establish equation  \eqref{eq:ninfty2}
in the distributional sense, we just pass to the limit, by weak-strong convergence, in equation \eqref{eq:n} .
On the other hand, using the definition of $p_k$ in \eqref{eq:p}, we have
$$
n_k p_k=\frac{k}{k-1}n_k^k =
\Big(1-\f{1}{k}\Big)^{1/(k-1)} p_k^{k/(k-1)} \underset{k \to \infty}{\longrightarrow}  p_\infty.
$$
Taking the limit $k\to \infty$, we deduce the monotone graph property
\beq\label{pninfty}
p_\infty(1-n_\infty)=0.
\eeq

In order to show the equivalence of \eqref{eq:ninfty} and \eqref{eq:ninfty2}, we need to prove that $\nabla p_{\infty}=n_\infty\nabla p_\infty$.
This es seen to be  equivalent to
$p_\infty \nabla n_\infty = 0$ by using the Leibnitz  rule in $H^1(\R^d)$ for \eqref{pninfty}. To prove the latter identity, we first write
$$
p_k\nabla n_k = \f{k}{k-1} n_k^k \nabla n_k = \f{\sqrt{k}}{k-1} n_k^{(k+1)/2}
\big(\sqrt{k}\,n_k^{(k-1)/2} \nabla n_k\big).
$$
From estimate \eqref{eq:bornpn}, the term between parentheses is uniformly
bounded in $L^2(Q_T)$ and since $(n_k)_k$ is uniformly (in $k$) bounded
in $L^\infty(Q_T)$, we conclude that
$$
\lim_{k\to \infty} \| p_k \nabla n_k\|_{L^2(Q_T)} = 0.
$$
We deduce then from the strong convergence of $(p_k)_k$ and the weak
convergence of $(\nabla n_k)_k$ that
\beq\label{pgradn}
p_\infty \nabla n_\infty = 0,
\eeq
as desired.

\medskip

\noindent\emph{Time continuity and initial trace. } Time continuity for the limit density $n_\infty$ follows from the monotonicity and the equation, as in the case $\nu=0$. Once we have continuity, the identification of the initial trace will follow from the equation for $n_k$, letting first $k\to\infty$ and then $t\to0$; see~\cite{PQV} for the details.

\medskip

\noindent\emph{Remark. } Since $p_\infty\ge0$,~\eqref{pgradn} implies that
\begin{equation}
\label{eq:grad.p.grad.n}
\nabla p_\infty\cdot \nabla n_\infty = 0.
\end{equation}
\section{The equation on $p_\infty$}
\label{sec:hs}
\setcounter{equation}{0}

In this section we give a rigorous derivation of equation \eqref{eqp2}, which is the most delicate point in the proof of~Theorem~\ref{th1}.

\medskip

\noindent(i) Our first goal is to establish that, in the weak sense,
\beq\label{ineq1b}
p_\infty \Delta p_\infty  +p_\infty G(p_\infty) \leq 0.
\eeq
Thanks to~\eqref{pgradn} and~\eqref{eq:grad.p.grad.n}, this is equivalent to proving that
\beq\label{ineq1}
p_\infty \Delta\big(p_\infty+\nu n_\infty\big)
+p_\infty G(p_\infty) \leq 0.
\eeq
In order to prove the latter inequality, we follow an idea of \cite{PQV} and use
a time regularization method {\it \`a la Steklov}. To this aim, we introduce
 a regularizing kernel $\omega_\eps(t)\geq 0$ with compact support
of length~$\eps$.

Let $n_{k,\eps}=n_k*\omega_\eps$.
From equation \eqref{eq:n}, we deduce
\beq\label{nkeps}
\p_t n_{k,\eps} - \Delta \omega_\eps*(n_k^k+\nu n_k) = (n_kG(p_k))*\omega_\eps.
\eeq
Then, for fixed $\eps>0$, $\Delta \omega_\eps*(n_k^k+\nu n_k)$
is bounded in $L^q(Q_T)$ for all $q\geq 1$.
Thus, we can extract a subsequence such that $(\nabla\omega_\eps*(n_k^k+\nu n_k))_k$
converges strongly in $L^2(Q_T)$.
Since we have strong convergence of $(n_k^k+\nu n_k)_k$ towards
$p_\infty+\nu n_\infty$, we deduce that the strong limit of
$(\nabla\omega_\eps*(n_k^k+\nu n_k))_k$ is equal to
$\nabla\omega_\eps*(p_\infty+\nu n_\infty)$.

Multiplying equation \eqref{nkeps} by $p_k$, we have
$$
p_k\p_tn_{k,\eps} = p_k \Delta\big(n_k^k*\omega_\eps+\nu n_{k,\eps}\big)
+p_k\big((n_kG(p_k))*\omega_\eps\big).
$$
We can pass to the limit $k\to\infty$ to get
$$
\lim_{k\to \infty} p_k\p_tn_{k,\eps} = p_\infty
\Delta\big(\omega_\eps*(p_\infty+\nu n_\infty)\big)
+p_\infty\big((n_\infty G(p_\infty))*\omega_\eps\big).
$$
To determine the sign, we decompose the left hand side term, divided by the harmless factor $k/(k-1)$,  as
$$
\begin{array}{l}
\displaystyle\int_\R n_k^{k-1}(t)\p_tn_k(s)\omega_\eps(t-s)\,ds =\\[4mm]
\qquad\qquad\displaystyle\underbrace{\int_\R n_k^{k-1}(s) \p_tn_k(s)\omega_\eps(t-s)\,ds}_{\mathcal{A}_k} + \underbrace{\int_\R (n_k^{k-1}(t)-n_k^{k-1}(s)) \p_tn_k(s) \omega_\eps(t-s)\,ds}_{\mathcal{B}_k}.
\end{array}
$$
On the one hand we have
$$
\mathcal{A}_k= \f 1k \int_\R
\p_tn^k(s) \omega_\eps(t-s)\,ds \to 0 \quad \mbox{ when } k\to \infty.
$$
As for  $\mathcal{B}_k$, we recall that $\p_tn_k\geq 0$ provided $\p_t n^{ini}\geq 0$; see Lemma~\ref{estim}. Thus, for $s>t$ we have $n_k^{k-1}(t)-n_k^{k-1}(s)\leq 0$.
Then, choosing $\omega_\eps$ such that $\mbox{supp }\omega_\eps\subset \R_-$,
we deduce that
$\mathcal{B}_k  \leq 0$,
which yields
$$
p_\infty \Delta\big(\omega_\eps*(p_\infty+\nu n_\infty)\big)
+p_\infty\big(n_\infty G(p_\infty)*\omega_\eps\big) \leq 0.
$$
It remains to pass to the limit $\eps\to 0$ in the regularization process.
We can pass to the limit in the weak formulation since we already
know that $\nabla p_\infty\in L^2(Q_T)$. Then, using
\eqref{pninfty}, we get the inequality \eqref{ineq1} and thus \eqref{ineq1b}.

\medskip

\noindent(ii)  Our second purpose is to establish the other inequality, namely
\beq\label{ineq2b}
p_\infty \Delta p_\infty
+p_\infty G(p_\infty) \geq  0 .
\eeq
To prove it, we multiply equation
\eqref{eq:pk} by a nonnegative test function $\phi(x,t)$  and integrate, and obtain
$$
\bea
\dis \iint_{Q_T} \phi \big(p_k\Delta p_k +p_k G(p_k) &- \nu \f{k-2}{k-1}
\frac{\nabla p_k\cdot\nabla n_k}{n_k}\big)
\\ [5pt]
& \dis = \f{1}{k-1}\iint_{Q_T} \left[
\phi\big(\partial_t p_k-|\nabla p_k|^2\big) +  \nu \nabla \phi\cdot \nabla p_k\right].
\eea
$$
From the proved bounds, the right hand side of the above equation
converges to $0$ as $k\to \infty$.
We can use integration by parts and rewrite the left hand side as
$$
\iint_{Q_T} \left(\phi p_k G(p_k)-p_k\nabla \phi\cdot \nabla p_k-\phi |\nabla p_k|^2
-\phi\nu \f{k(k-2)}{k-1} n_k^{k-3}|\nabla n_k|^2\right).
$$
Since the last term is nonpositive, we obtain that
$$
\liminf_{k\to\infty}\iint_{Q_T} \left(\phi p_k G(p_k)-p_k\nabla \phi\cdot \nabla p_k-\phi |\nabla p_k|^2\right) \geq 0.
$$
From weak-strong convergence in products, or convexity inequalities in the weak limit, we finally conclude
$$
\iint_{Q_T}\left(\phi p_\infty G(p_\infty) - p_\infty\nabla\phi\cdot \nabla p_\infty
-\phi |\nabla p_\infty|^2 \right)\geq 0.
$$
This is the weak formulation of  \eqref{ineq2b}.

\medskip

\noindent\emph{Remark. } A careful inspection of the proof of~\eqref{ineq2b} shows that~\eqref{eqp1} holds if and only if  $\nabla p_k$ converges strongly in $L^2(Q_T)$ and $k n_k^{k-3}|\nabla n_k|^2$ converges weakly to 0 locally in $L^1(Q)$. Since we have proved~\eqref{eqp1}, we conclude that we have the two mentioned convergence results.

\section{Uniqueness for the limit model}
\label{sec:uniq}
\setcounter{equation}{0}

In this section we prove that the limit problem \eqref{eq:ninfty}
admits at most one solution. We will adapt Hilbert's duality method
in the spirit of \cite{PQV}.
\begin{theorem}\label{th:uniq}
Let $T>0$, $\nu>0$.
There  is  a unique pair $(n,p)$ of functions
in $L^\infty([0,T];L^1(\mathbb{R}^d)\cap L^\infty(\R^d))$, $n\in C([0,T]; L^1(\mathbb{R}^d))$, $n(0)=n^{\text{ini}}$, $p\in P_\infty(n)$,  satisfying~\eqref{eq:ninfty} in the sense of distributions and such that
$\nabla n,\;\nabla p \in L^2(Q_T)$,
$\p_t n,\;\p_t p \in {\cal M}^1(Q_T)$.
\end{theorem}

\noindent\emph{Proof. }
Let us consider two solutions $(n_1,p_1)$ and $(n_2,p_2)$. Then
for any test function $\phi$ with $\phi\in W^{2,2}(Q_T)$
and $\p_t\phi \in L^2(Q_T)$, we have
\begin{equation}\label{eq:uniq1}
\iint_{Q_T} \Big( (n_1-n_2)\p_t\phi + (p_1-p_2+\nu (n_1-n_2)) \Delta\phi
+ \big(n_1G(p_1)-n_2G(p_2)\big) \phi \Big)= 0,
\end{equation}
which can be rewritten as
\begin{equation}\label{eq:uniq2}
\iint_{Q_T} \big(\nu(n_1-n_2)+p_1-p_2\big)
\big( A\p_t \phi + \Delta \phi + A G(p_1) \phi - B \phi \big) = 0,
\end{equation}
where
$$
\begin{array}{c}
\dis 0 \leq A = \f{n_1-n_2}{\nu(n_1-n_2)+p_1-p_2} \leq \f{1}{\nu},  \\[3mm]
\dis 0 \leq B = -n_2\f{G(p_1)-G(p_2)}{\nu(n_1-n_2)+p_1-p_2} \leq \kappa,
\end{array}
$$
for some nonnegative constant $\kappa$. To arrive to these bounds on $A$
we set $A=0$ when $n_1=n_2$, even if $p_1=p_2$.
Since $A$ can vanish, we use a smoothing argument by introducing the
regularizing sequences $(A_n)_n$, $(B_n)_n$ and $(G_{1,n})_n$ such that
$$
\begin{array}{lll}
\dis \|A-A_n\|_{L^2(Q_T)} < \alpha/n, \qquad & \dis 1/n < A_n \leq 1, \\[2mm]
\dis \|B-B_n\|_{L^2(Q_T)} < \beta/n, \qquad & \dis 0 \leq B_n \leq \beta_2,
\qquad \|\p_t B_n\|_{L^1(Q_T)} \leq \beta_3,  \\[2mm]
\dis \|G_{1,n} -G(p_1)\|_{L^2(Q_T)} \leq \delta/n, \qquad & \dis |G_{1,n}| < \delta_2,
\qquad \| \nabla G_{1,n} \|_{L^2(Q_T)} \leq \delta_3,
\end{array}
$$
for some nonnegative constants $\alpha$, $\beta$, $\beta_2$,
$\beta_3$, $\delta$, $\delta_2$, $\delta_3$.

Given any arbitrary smooth function $\psi$ compactly
supported, we consider the solution $\phi_n$ of the backward heat equation
\begin{equation}\label{eq:phin}
\left\{\begin{array}{l}
\dis \p_t \phi_n + \f{1}{A_n} \Delta \phi_n + G_{1,n} \phi_n -
\f{B_n}{A_n} \phi_n = \psi\qquad \mbox{in } Q_T, \\[3mm]
\dis \phi_n(T) = 0.
\end{array}\right.
\end{equation}
The coefficient $1/A_n$ is continuous, positive and bounded below away from
zero. Then the equation satisfied by $\phi_n$ is parabolic. Hence
$\phi_n$ is smooth and since $\psi$ is compactly supported, we
have that $\phi_n$, $\Delta \phi_n$ and therefore $\p_t\phi_n$
are $L^2$-integrable.
Therefore, we can use $\phi_n$ as a test function in \eqref{eq:uniq2}.
Then, by  the definition of $A$, we have
$$
\iint_{Q_T} (n_1-n_2)\psi = \iint_{Q_T}
\big(\nu(n_1-n_2)+p_1-p_2\big) A \psi.
$$
Inserting \eqref{eq:phin} and substracting \eqref{eq:uniq2}, we obtain
$$
\iint_{Q_T} (n_1-n_2)\psi  = I_{1n} + I_{2n} + I_{3n},
$$
where
$$
\begin{array}{l}
\dis I_{1n} = \iint_{Q_T}
\big(\nu(n_1-n_2)+p_1-p_2\big)\Big( \big(\f{A}{A_n}-1\big)
\big(\Delta \phi_n- C_n\phi_n\big)\Big),  \\[3mm]
\dis I_{2n} = \iint_{Q_T} \big(\nu(n_1-n_2)+p_1-p_2\big)
(B-B_n)\phi_n,   \\[3mm]
\dis I_{3n} = \iint_{Q_T} (n_1-n_2)
\big(G_{1,n} - G(p_1)\big) \phi_n.
\end{array}
$$

The convergence towards  $0$ of the terms
$I_{in}$, $i=1,2,3$ is now a consequence on some estimates on the test functions $\phi_n$ which are gathered in Lemma~\ref{lem:techniq} below. Indeed, applying the mentioned estimates and  Cauchy-Schwarz inequality we have
$$
\begin{array}{l}
I_{1n} \leq K \|(A-A_n)/\sqrt{A_n}\|_{L^2(Q_T)}
\leq K \sqrt{n} \|A-A_n\|_{L^2(Q_T)} \leq K \alpha/\sqrt{n},\\
I_{2n} \leq K \|B-B_n\|_{L^2(Q_T)} \leq K\gamma/n,\\
I_{3n} \leq K \delta /n,
\end{array}
$$
(in all the computations, $K$ denotes various nonnegative constants).
Then letting $n\to \infty$, we conclude that
$$
\iint_{Q_T} (n_1-n_2)\psi = 0,
$$
for any smooth function $\psi$ compactly supported, hence $n_1=n_2$.
It is then obvious, thanks to \eqref{eq:uniq1}, that $p_1=p_2$.
\qed

\begin{lemma}\label{lem:techniq}
Under the assumptions of Theorem \ref{th:uniq}, we have the uniform
bounds, only depending on $T$ and $\psi$,
$$
\|\phi_n\|_{L^\infty(Q_T)} \leq \kappa_1, \quad
\sup_{0\leq t\leq T} \|\nabla \phi_n(t)\|_{L^2(\R^d)} \leq \kappa_2, \quad
\|1/\sqrt{A_n} (\Delta \phi_n -B_n\phi_n) \|_{L^2(Q_T)} \leq \kappa_3.
$$
\end{lemma}
\noindent\emph{Proof. }
The first bound is a consequence of the maximum principle on \eqref{eq:phin}.
Then multiplying \eqref{eq:phin} by $\Delta \phi_n-B_n\phi_n$ and integrating
on $\R^d$, we get
$$
\begin{array}{l}
\dis - \f{1}{2} \f{d}{dt} \int_{\R^d} |\nabla \phi_n(t)|^2
- \f{1}{2} \f{d}{dt} \int_{\R^d} B_n \phi_n^2(t) + \int_{\R^d} \f{1}{A_n}
|\Delta \phi_n -B_n\phi_n|^2(t) + \f{1}{2} \int_{R^d} (\p_tB_n \phi_n^2)(t),
\\[3mm]
\dis = \int_{\R^d} \Big( G_{1,n} |\nabla \phi_n|^2 +
 \phi_n \nabla \phi_n\cdot \nabla G_{1,n} +
B_n G_{1,n} \phi_n^2 + (\Delta \psi -B_n\psi) \phi_n \Big)(t).
\end{array}
$$
After an integration in time on $[t,T]$, we deduce
$$
\f{1}{2}\|\nabla \phi_n(t)\|_{L^2(\R^d)} + \int_t^T \int_{\R^d} \f{1}{A_n}
|\Delta \phi_n -B_n\phi_n|^2 \leq K \Big( 1-t+ \int_t^T
\|\nabla \phi_n(s)\|_{L^2(\R^d)}\,ds\Big),
$$
where we use the bounds on $\nabla G_{1,n}$ and $\p_tB_n$ by construction
of the regularization.
We conclude by applying Gronwall's Lemma.
\qed

\section{Further regularity and velocity of the free boundary}
\label{sec:fb.velocity}
\setcounter{equation}{0}

Remember that both $p_\infty$ and $n_\infty$ belong to $H^1(\mathbb{R}^d)$ for almost every $t>0$. This regularity cannot be improved, because there are jumps in the gradients of   both $p_\infty$ and $n_\infty$ at the free boundary. As a consequence,  their laplacians are not functions, but measures.  However, these singularities cancel in the combination $\Sigma_\infty = p_\infty+\nu n_\infty$, as we will see now.
\begin{lemma} With the assumptions of Theorem \ref{th1}, the quantity $\Sigma_\infty$ belongs to $ L^2((0,T); H^2(\R^d))$ for all $T >0$ and we have the estimate
$$
\iint_{Q_T} (\Delta\Sigma_\infty)^2  \le C(T).
$$
\end{lemma}
\noindent\emph{Proof. }
We recall the definition of $\Sigma_k$ in \eqref{def:sigmak}. Since  $\nabla\Sigma_k=n_k\nabla p_k+ \nu\nabla n_k$, estimate~\eqref{eq:bornpn} yields that for all $0< t \leq T$,
$$
\int_{\R^d}|\nabla \Sigma_k(t)|^2 \le C(T).
$$
We now multiply the equation~\eqref{eq:sigma} by $\Delta \Sigma_k$, and integrate in $Q_T$, $0< T<\infty$, to obtain, using that $\Sigma_k'>\nu$ and the fact that both $n_k$ and $G(p_k)$ are nonnegative,
$$
 \iint_{Q_T}(\Delta\Sigma_k)^2\le\frac12 \int_{\mathbb{R}^d}|\nabla \Sigma_k|^2(0) +C(T).
$$
The result follows directly.
\qed

This implies in particular that in the limit
$\Sigma_\infty(\cdot,t)\in H^2(\mathbb{R}^d)$ for almost every $t>0$. Hence, the size of the jump (downwards) of
$\nabla p_\infty$ at the free boundary coincides with  the size of the  jump (upwards) of $\nu\nabla n_\infty$ there.

Concerning the time regularity, the limit equation for the density \eqref{eq:sigmak}, now tells us that $\partial_t n_\infty\in L^2(Q_T)$. Hence $n_\infty \in H^1(Q_T)$. We do not have a similar property for the pressure (think of the situation when two tumors meet).

\bigskip


Our last goal is to derive formally an asymptotic value for the free boundary speed in a particular example. Let
 $\Omega(t)$ denote, as before, the space filled by the tumor at time $t$.
We notice that $n_\infty$ solves
$$
\partial_t n_\infty=\nu\Delta n_\infty+G(0)n_\infty,\qquad x\in\mathbb{R}^d\setminus\Omega(t),\; t>0,
$$
with boundary conditions
$$
n_\infty=1, \quad \nu\partial_n n_\infty=\partial_n p_\infty, \qquad x\in\partial\Omega (t),\; t>0.
$$
If $\Omega(t)$ were known, the problem would be overdetermined. This is precisely what fixes the dynamics of the free boundary. Let us assume that the tumor is  a ball centered at the origin,
$$
\Omega(t)=\{x: p_\infty(x,t)>0\} = \{x: n_\infty(x,t)=1\} = B_{R(t)}(0).
$$
We look for a solution which is spherically symmetric
$n_\infty(r,t)$, $p_\infty(r,t)$.
We set $\sigma=R'(t)$. In opposition to other models of tumor growth (see \cite{TVCVDP} for instance), here there are no  radial solutions with constant speed. However, following \cite{PQV} Appendix A, we expect  our solution to behave for large times as a one dimensional traveling wave (with constant speed).

In order to analyze the expected asymptotic constant speed, we set $n_R(r-\sigma t)=n_\infty(r,t)$ and $p_R(r-\sigma t)=p_\infty(r,t)$.
Introducing this ansatz in  equation \eqref{eq:ninfty}, we obtain
\begin{equation}
\label{eq:tw}
-\sigma n_R' =  p_R''+ \frac{d-1}{r}p_R' + \nu n_R''
+\nu \frac{d-1}{r} n_R' + n_R G(p_R).
\end{equation}
On $\R^d \setminus \Omega(t)$, we have $p_\infty=0$, then
integrating~\eqref{eq:tw} in $(R(0),\infty)$, we get
$$
\sigma n_R(R(0)) = -\nu n_R'(R(0)^+)+ \nu (d-1) \int_{R(0)}^\infty \frac{n_R'}{r}dr
+ G(0)\int_{R(0)}^\infty n_R dr.
$$
In a one dimensional setting ($d=1$)
and using the boundary relation at the interface of $\Omega(0)$, we deduce
\beq\label{eqTWsigma}
\sigma = - p_R'(R(0)^-) + G(0) \int_{R(0)}^{\infty} n_R(r) dr.
\eeq
We recall that for the classical Hele-Shaw model without active motion
(i.e. $\nu=0$), the traveling velocity is $\sigma_0=-p_R'(R(0)^-)$.
Since $n_R(R(0))=1$ and $n_R$ is continuous and nonnegative, we have
$\int_{R(0)}^{\infty} n_R(r)dr >0$. Then
we conclude from equation \eqref{eqTWsigma} that $\sigma>\sigma_0$.

\begin{figure}
\begin{center}
\includegraphics[width=0.47\textwidth]{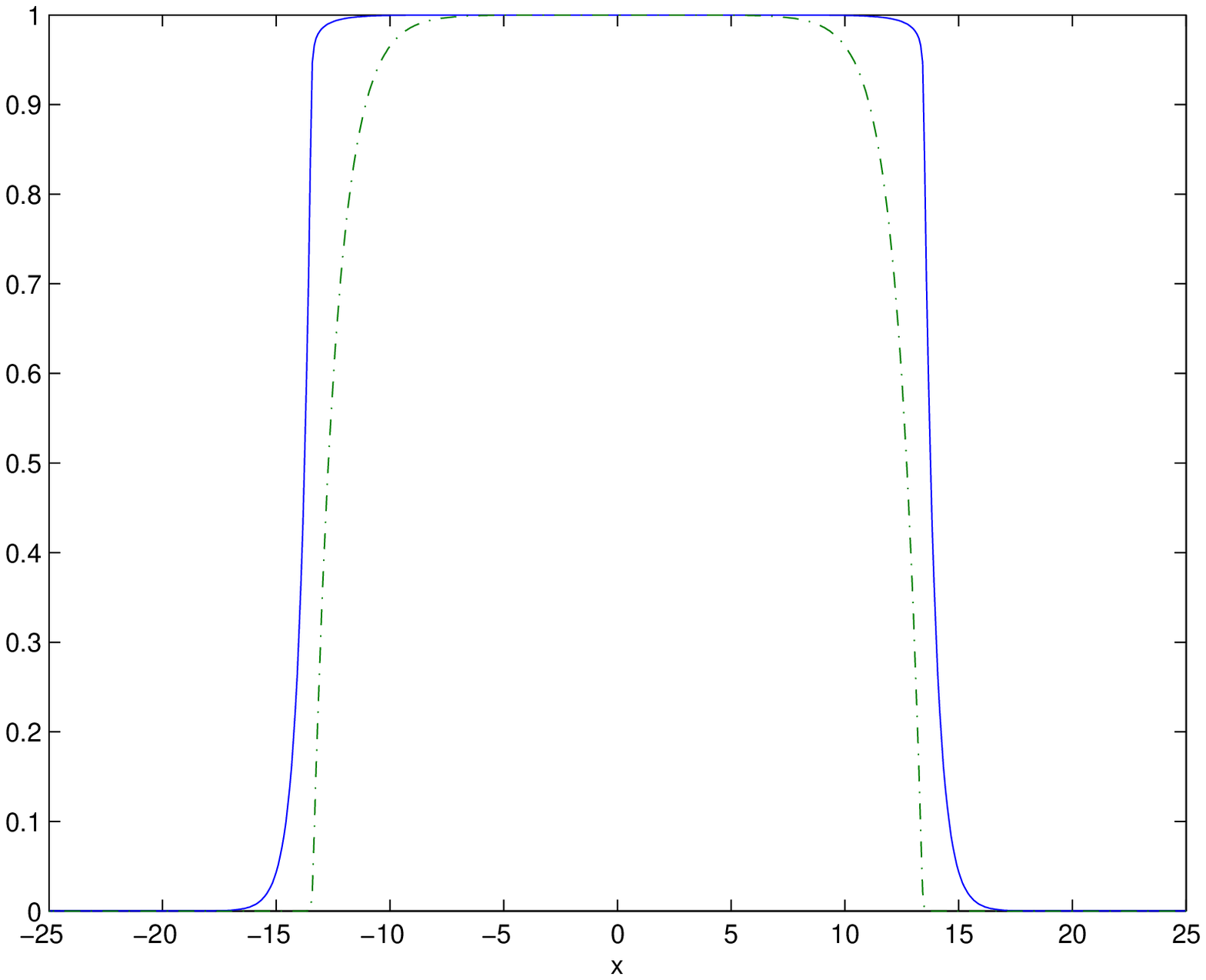}
\includegraphics[width=0.47\textwidth]{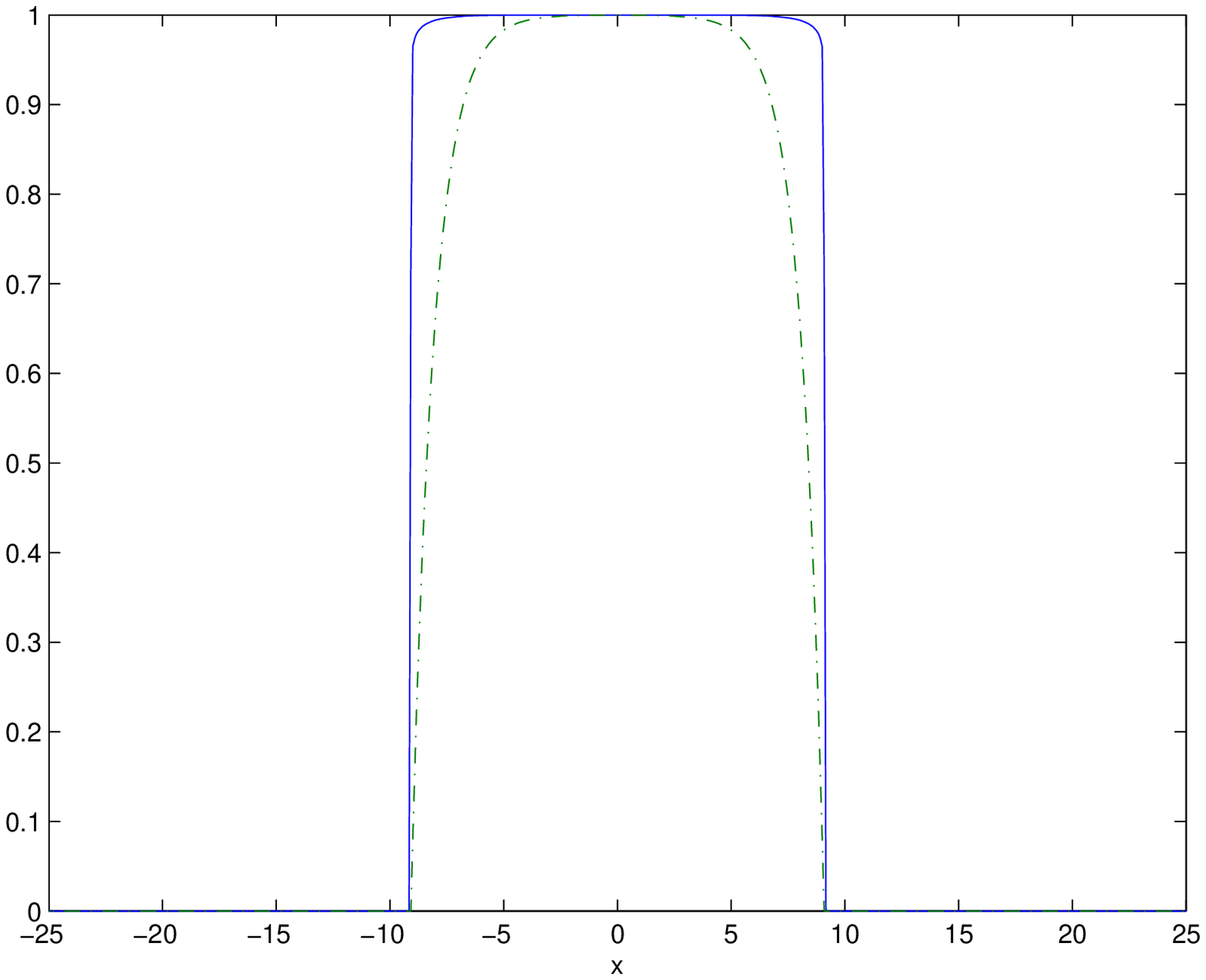}
\vspace{-7mm}
\caption{Shape of the traveling waves obtained thanks to a numerical
discretization of the system \eqref{eq:n}--\eqref{eq:p} with $k=100$
and $\nu=0.5$ (left) or $\nu=0$ (right) for the same initial data and the same final time.
The density $n$ is plotted in line whereas the pressure is represented in
dashed line.
We notice the regularity of $n$ in the case $\nu=0.5$, whereas it has
a jump at the interface when  $\nu=0$. Also the free boundary moves faster when active motion is present.
}\label{fig:TWact}
\end{center}
\end{figure}

We can do a more precise computation confirming the above statement for the one-dimensional case. From the complementarity relation \eqref{eqp2},
we have $-p_R'' = G(p_R)$ on $\Omega(0)$.
Multiplying this latter equation by $p_R'$
and integrating on $(0,R(0))$, we deduce
$$
p_R'(R(0)^-)^2 = 2 \int_0^{R(0)} p_R' G(p_R) dr.
$$
In the center of the tumor, we expect a maximal packing of the cells. Therefore, we have the boundary conditions
$$
\lim _{r\to 0} p_R(r) = P_M, \quad \lim_{r\to 0} p_R'(r)=0.
$$
Since $p_R'' = -G(p_R)\leq 0$, we deduce that $p_R' < 0$ and we can make the
change of variable
$$
p_R'(R(0)^-)^2 = 2 \int_0^{R(0)} p_R' G(p_R)\,dr = 2\int_0^{P_M} G(q)\,dq.
$$
The quantity $\sigma_0=\sqrt{2\int_0^{P_M} G(q)\,dq}$ is the
traveling velocity for a tumor spheroid in the case $\nu=0$; see Appendix A.1
of \cite{PQV}.
Combining this with \eqref{eqTWsigma}, we deduce that the growth of the tumor is faster
with active motion than in the case $\nu=0$.

In Figure \ref{fig:TWact}, we display numerical simulations obtained
from a discretization with a finite volume scheme of system \eqref{eq:n}--\eqref{eq:p}
for $k=100$. The left picture presents the result for $\nu=0.5$, and
the right for $\nu=0$ (i.e. without active motion). We use the growth
function $G(p)=1-p$ and the results in both cases with the same initial data
and at final time $t=10$. We notice that in the case $\nu=0.5$ the density function
is smooth and the domain occupied by the tumor is larger than in the
case without active motion, which suggests as explained above a faster
invasion speed.

%
%
%


\begin{thebibliography}{99}

\bibitem{Bellomo1} Bellomo, N.; Li, N.~K.; Maini, P.~K.  \emph{On
the foundations of cancer modelling: selected topics, speculations,
and perspectives}. Math. Models Methods Appl. Sci. 18  (2008),
no.~4, 593--646.

\bibitem{BI2} B\'enilan, Ph.; Igbida, N. \emph{La limite de la solution de
$u_t=\Delta_p u^m$ lorsque $m\to\infty$}. {C.\ R.\ Acad.\ Sci.\
Paris S\'er.\ I Math.} 321 (1995), no.~10, 1323--1328.



\bibitem{BOBAM}  Betteridge, R.; Owen, M.~R.; Byrne, H.~M.; Alarc\'on, T.; Maini, P.~K. \emph{The impact of cell crowding and active cell movement on vascular tumour growth}. Netw. Heterog. Media 1 (2006), no.~4, 515--535.

\bibitem{Bru} Br\'u, A.; Albertos, S.;  Subiza, J.~L.;  Asenjo,  J.~A.;
Br{\oe}, I. \emph{The universal dynamics of tumor growth}. Biophys. J.
85  (2003), no.~5, 2948--2961.

\bibitem{byrne-chaplain} Byrne, H.~M.; Chaplain, M.~A. \emph{Growth
of necrotic tumors in the presence and absence of inhibitors}. Math.
Biosci. 135 (1996), no.~15, 187--216.

\bibitem{byrne-drasdo} Byrne, H.~M.; Drasdo, D. \emph{Individual-based
and continuum models of growing cell populations: a comparison}. J.
Math. Biol. 58 (2009), no. 4-5, 657--687.

\bibitem{byrne-preziosi} Byrne, H.~M.; Preziosi, L. \emph{Modelling
solid tumour growth using the theory of mixtures}. Math. Med. Biol.
20 (2003), no.~4, 341--366.

\bibitem{chaplain} Chaplain, M.~A. J. \emph{Avascular growth,
angiogenesis and vascular growth in solid tumours: the mathematical
modeling  of the stages of tumor development}. Math. Comput.
Modeling 23 (1996), no.~6, 47--87.


\bibitem{ciarletta} Ciarletta, P.; Foret, L.; Ben Amar, M. \emph{The radial growth phase of malignant melanoma: multi-phase modelling, numerical simulations and linear stability analysis}. J. R. Soc. Interface 8 (2011) no.~56, 345--368.


\bibitem{cui} Cui, S. \emph{Formation of necrotic cores in the growth of tumors: analytic results}.
Acta Math. Sci. Ser. B Engl. Ed. (2006), no.~4, 781--796.

\bibitem{cui_escher} Cui, S.; Escher, J. \emph{Asymptotic behaviour of
solutions of a multidimensional moving boundary problem modeling
tumor growth}. Comm. Partial Differential Equations 33 (2008),
no.~4--6, 636--655.

\bibitem{Drasdo_H} Drasdo, D.;  Hoehme, S.  \emph{Modeling the impact of granular embedding media,
and pulling versus pushing cells on growing cell
clones}. New J. Phys. 14 (2012) 055025 (37pp).


\bibitem{friedman_hu} Friedman, A.; Hu, B.  \emph{Stability and instability of
Liapunov-Schmidt and Hopf bifurcation for a free boundary problem
arising in a tumor model}. Trans. Am. Math. Soc. 360 (2008), no.~10,
5291--5342.

\bibitem{GQ2} Gil, O.; Quir\'os, F. \emph{Boundary layer formation in the transition from the porous
media equation to a Hele-Shaw flow}. Ann. Inst. H. Poincar\'{e}
Anal. Non Lin\'{e}aire 20 (2003), no. 1, 13--36.


\bibitem{greenspan} Greenspan, H.~P. \emph{Models for the growth of a solid tumor by diffusion}.
Stud. Appl. Math. 51  (1972), no.~4, 317--340.


\bibitem{Lowengrub_survey}   Lowengrub, J.~S.;  Frieboes H.~B.;  Jin, F.;  Chuang, Y.-L.;  Li, X.;
 Macklin, P.;  Wise, S.~M.;  Cristini, V.  \emph{Nonlinear modelling of cancer: bridging the gap
between cells and tumours}. Nonlinearity 23 (2010), no.~1, R1--R91.


\bibitem{PQV} Perthame, B.; Quir\'os, F.;  V\'azquez, J.~L. \emph{The Hele-Shaw asymptotics for mechanical models of tumor growth}. Arch. Ration. Mech. Anal., to appear.

\bibitem{preziosi_tosin}  Preziosi, L.;  Tosin, A. \emph{Multiphase modelling of tumour growth
and extracellular matrix interaction:  mathematical tools and
applications}. J. Math. Biol. 58 (2009),  no.~4-5, 625--656.

\bibitem{JJP} Ranft, J.; Basana, M.;  Elgeti, J.;  Joanny, J.-F.;  Prost, J.; J\"ulicher, F.
\emph{Fluidization of tissues by cell division and apoptosis}. Proc.
Natl. Acad. Sci. USA (2010), no.~49, 20863--20868.


\bibitem{RCM}  Roose, T.; Chapman, S.~J.; Maini, P.~K. \emph{Mathematical models of avascular tumor growth}. SIAM Rev. 49 (2007), no.~2, 179--208.

\bibitem{SLCF} Saut, O.; Lagaert, J.-B.; Colin, T.; Fathallah-Shaykh, H.~M. \emph{ A multilayer grow-or-go model for GBM: effects of
invasive cells and anti-angiogenesis on growth}. Preprint 2012.

\bibitem {TVCVDP} Tang, M.; Vauchelet, N.; Cheddadi, I.; Vignon-Clementel, I.; Drasdo, D.;
Perthame, B. \emph{Composite waves for a cell population system
modelling tumor growth and invasion}. Chin. Ann. Math.  Ser. B 34
(2013), no.~2, 295--318.

\bibitem{JLV} V\'{a}zquez, J.~L. \lq\lq The porous medium
equation. Mathematical theory''. Oxford Mathematical Monographs. The
Clarendon Press, Oxford University Press, Oxford, 2007. ISBN:
978-0-19-856903-9.

\end{thebibliography}
\end{document}